\documentclass[12pt]{amsart}
\usepackage{amsxtra}
\usepackage{amssymb, amsmath}
\usepackage{amsthm}
\usepackage {hyperref}

\newtheorem{theorem}{Theorem}
\newtheorem{remark}{Remark}
\newtheorem{example}{Example}

\thanks{This research is partially supported by the Natural Sciences and
Engineering Research Council's Discovery Grant to Pranesh Kumar}


\begin{document}
\title[Information Inequalities]{Generalized Non-Symmetric Divergence
Measures and Inequalities}

\author{Inder Jeet Taneja}
\address{Inder Jeet Taneja\\
Departamento de Matem\'{a}tica\\
Universidade Federal de Santa Catarina\\
88.040-900 Florian\'{o}polis, SC, Brazil}
\email{taneja@mtm.ufsc.br} \urladdr{http://www.mtm.ufsc.br/$\sim
$taneja}

\author{Pranesh Kumar}
\address{Pranesh Kumar\\
Mathematics Department\\
College of Science and Management\\
University of Northern British Columbia\\
Prince George BC V2N4Z9, Canada.} \email{kumarp@unbc.ca}
\urladdr{http://web.unbc.ca/$\sim$kumarp}

\begin{abstract}
In this paper we consider one parameter generalizations of some
non - symmetric divergence measures. Measures are \textit{relative
information}, $\chi ^2 - $\textit{divergence}, \textit{relative
J-divergence}, \textit{relative Jensen-Shannon divergence }and
\textit{relative arithmetic and geometric divergence}. All the
generalizations considered can be written as particular cases of
Csisz\'{a}r \textit{f-divergence}. By conditioning the probability
distributions, relationships among the \textit{relative divergence
measures }are obtained.

\end{abstract}

\maketitle

\section{Introduction}

Let
\[
\Gamma _n = \left\{ {P = (p_1 ,p_2 ,...,p_n )\left| {p_i > 0,\sum\limits_{i
= 1}^n {p_i = 1} } \right.} \right\},
\quad
n \geqslant 2,
\]
\noindent be the set of all complete finite discrete probability
distributions. There are many information and divergence measures
given in the literature on information theory and statistics. Some
of these are symmetric with respect to probability distribution,
while others are not. In this paper, we work only with
non-symmetric measures. Throughout the paper it is under stood
that the probability distributions $P,Q \in \Gamma _n $.

\subsection*{1.1. Non-Symmetric Divergence Measures }

The followings are some non-symmetric measures of information, the
most famous among them being $\chi ^2 - $\textit{divergence} and
Kullback-Leibler \textit{relative information}. All of the
following measures can be written in pairs by interchanging $P$
and $Q$, and $p_i $ and $q_i $. The measures are non-symmetric in
the sense that the expression changes when this interchange is
made.\\

\textbf{$\bullet$ $\chi ^2 - $Divergence} (Pearson \cite{pea})
\begin{equation}
\label{eq1} \chi ^2(P\vert \vert Q) = \sum\limits_{i = 1}^n
{\frac{(p_i - q_i )^2}{q_i }} = \sum\limits_{i = 1}^n {\frac{p_i^2
}{q_i } - 1}.
\end{equation}

\textbf{$\bullet$ Relative Information} (Kullback and Leiber
\cite{kul})
\begin{equation}
\label{eq2} K(P\vert \vert Q) = \sum\limits_{i = 1}^n {p_i \ln
(\frac{p_i }{q_i })}.
\end{equation}

\textbf{$\bullet$ Relative Jensen-Shannon Divergence} (Sibson
\cite{sib}, Sgarro \cite{sga})
\begin{equation}
\label{eq3} F(P\vert \vert Q) = \sum\limits_{i = 1}^n {p_i \ln
\left( {\frac{2p_i }{p_i + q_i }} \right)}.
\end{equation}

\textbf{$\bullet$ Relative Arithmetic-Geometric Divergence}
(Taneja \cite{tan7})
\begin{equation}
\label{eq4}
G(P\vert \vert Q) = \sum\limits_{i = 1}^n {\left( {\frac{p_i + q_i }{2}}
\right)\ln \left( {\frac{p_i + q_i }{2p_i }} \right)}
\end{equation}

\textbf{$\bullet$ Relative J-Divergence} (Dragomir et al.
\cite{dgp1})
\begin{equation}
\label{eq5} D(P\vert \vert Q) = \sum\limits_{i = 1}^n {(p_i - q_i
)\ln \left( {\frac{p_i + q_i }{2q_i }} \right)}.
\end{equation}

Symmetric versions of above measures are given by
\begin{equation}
\label{eq6}
\Psi (P\vert \vert Q) = \chi ^2(P\vert \vert Q) + \chi ^2(Q\vert \vert P),
\end{equation}
\begin{align}
\label{eq7} J(P\vert \vert Q) & = K(P\vert \vert Q) + K(Q\vert
\vert P)\\
& = D(P\vert \vert Q) + D(Q\vert \vert P),\notag
\end{align}
\begin{equation}
\label{eq8}
I(P\vert \vert Q) = \frac{1}{2}\left[ {F(P\vert \vert Q) + F(Q\vert \vert
P)} \right]
\end{equation}

\noindent and
\begin{equation}
\label{eq9} T(P\vert \vert Q) = \frac{1}{2}\left[ {G(P\vert \vert
Q) + G(P\vert \vert Q)} \right].
\end{equation}

After simplification, we can write
\begin{equation}
\label{eq10}
J(P\vert \vert Q) = 4\left[ {I(P\vert \vert Q) + T(Q\vert \vert P)} \right]
\end{equation}

\noindent and
\begin{equation}
\label{eq11}
D(Q\vert \vert P) = \frac{1}{2}\left[ {F(P\vert \vert Q) + G(P\vert \vert
Q)} \right].
\end{equation}

Dragomir et al. \cite{dsb} studied the measures (\ref{eq5}),
referred to subsequently as \cite{tan6} \textit{symmetric
chi-square divergence}. Measure (\ref{eq7}) is known as
Jeffreys-Kullback-Leiber \cite{jef, kul} \textit{J-divergence}.
The measure (\ref{eq8}) is \textit{Jensen-Shannon divergence}
studied by Sibson \cite{sib} and Burbea and Rao \cite{bur1,
bur2}). Measure (\ref{eq9}) is \textit{arithmetic and geometric
mean divergence }studied by Taneja \cite{tan2}. More details on
some of these measures can be found in Taneja \cite{tan1, tan2}
and in the on line book by Taneja \cite{tan4}.

In this paper our aim is to work with one parameter
generalizations of the \textit{non symmetric divergence measures}
given by (\ref{eq1})-(\ref{eq5}). We call these generalizations,
\textit{non-symmetric divergence measures of type s}. Also, we
call the measure $K(Q \vert \vert P)$ the \textit{adjoint} of $K(P
\vert \vert Q)$ and \textit{vice-versa}. The same is with the
other measures.

\section{Non-Symmetric Divergence Measures of Type s}
In this section we introduce one parameter generalizations of the
measures given by (\ref{eq1})-(\ref{eq5}). Generalization of the
measures (\ref{eq1}) and (\ref{eq2}) is already known in the
literature and has been studied by many authors. Here we refer to
it as \textit{relative information of type s}.\\

\textbf{$\bullet$ Relative Information of Type s}
\begin{equation}
\label{eq12} \Phi _s (P\vert \vert Q) = \begin{cases}
 {K_s (P\vert \vert Q) = \left[ {s(s - 1)} \right]^{ - 1}\left[
{\sum\limits_{i = 1}^n {p_i^s q_i^{1 - s} - 1} } \right],} & {s \ne 0,1} \\
 {K(Q\vert \vert P) = \sum\limits_{i = 1}^n {q_i \ln \left( {\frac{q_i }{p_i
}} \right)} ,} & {s = 0} \\
 {K(P\vert \vert Q) = \sum\limits_{i = 1}^n {p_i \ln \left( {\frac{p_i }{q_i
}} \right)} ,} & {s = 1} \\
\end{cases},
\end{equation}

\noindent for all $s \in \mathbb{R}$.

The measure (\ref{eq12}) admits the following particular cases:
\begin{itemize}
\item[(i)] $\Phi _{ - 1} (P\vert \vert Q) = \frac{1}{2}\chi
^2(Q\vert \vert P)$.
\item[(ii)] $ \Phi _0 (P\vert \vert Q) =
K(Q\vert \vert P)$.

\item[(iii)] $ \Phi _{1 / 2} (P\vert \vert Q) = 4\left[ {1 -
B(P\vert \vert Q)} \right] = 4h(P\vert \vert Q)$.

\item[(iv)] $\Phi _1 (P\vert \vert Q) = K(P\vert \vert Q)$.

\item[(v)] $\Phi _2 (P\vert \vert Q) = \frac{1}{2}\chi ^2(P\vert
\vert Q)$.
\end{itemize}

The measures $B(P\vert \vert Q)$ and $h(P\vert \vert Q)$ appearing in part
(iii) are given by
\begin{equation}
\label{eq13}
B(P\vert \vert Q) = \sqrt {p_i q_i } ,
\end{equation}

\noindent and
\begin{equation}
\label{eq14}
h(P\vert \vert Q) = 1 - B(P\vert \vert Q) = \frac{1}{2}\sum\limits_{i = 1}^n
{(\sqrt {p_i } - \sqrt {q_i } )^2}
\end{equation}

\noindent respectively. $B(P\vert \vert Q)$ is known as the
Bhattacharyya \cite{bha} coefficient and $h(P\vert \vert Q)$ as
Hellinger \cite{hel} discrimination.

From (\ref{eq12}), we observe that $\Phi _2 (P\vert \vert Q) =
\Phi _{ - 1} (Q\vert \vert P)$ and $\Phi _1 (P\vert \vert Q) =
\Phi _0 (Q\vert \vert P)$.\\

We now present new one parameter generalizations of the measures
given by (\ref{eq3})-(\ref{eq5}).\\

\textbf{$\bullet$ Unified Relative JS and AG -- Divergence of Type
s}

\bigskip
We consider the following unified one parameter generalization of
measures (\ref{eq3}) and (\ref{eq4}) simultaneously.
\begin{equation}
\label{eq15} \Omega _s (P\vert \vert Q) = \begin{cases}
 {FG_s (P\vert \vert Q) = \left[ {s(s - 1)} \right]^{ - 1}\left[
{\sum\limits_{i = 1}^n {p_i \left( {\frac{p_i + q_i }{2p_i }} \right)^s - 1}
} \right],} & {s \ne 0,1} \\
 {F(P\vert \vert Q) = \sum\limits_{i = 1}^n {p_i \ln \left( {\frac{2p_i
}{p_i + q_i }} \right)} ,} & {s = 0} \\
 {G(P\vert \vert Q) = \sum\limits_{i = 1}^n {\left( {\frac{p_i + q_i }{2}}
\right)\ln \left( {\frac{p_i + q_i }{2p_i }} \right)} ,} & {s = 1} \\
\end{cases}.
\end{equation}

The adjoint of $\Omega _s (P\vert \vert Q)$ written as $\Omega _s
(Q\vert \vert P)$ is obtained by interchanging $P$ and $Q$, and
$p_i $ and $q_i $ in the expression (\ref{eq15}). The measures
$\Omega _s (Q\vert \vert P)$ can also be obtained from
(\ref{eq12}) by replacing $p_i $ by $\frac{p_i + q_i }{2}$.\\

We have the following particular cases of $\Omega _s (P\vert \vert
Q)$ and $\Omega _s (Q\vert \vert P)$:

\begin{itemize}\item[(i)]
$\Omega _{ - 1} (P\vert \vert Q) = \Omega _{ - 1} (Q\vert \vert P)
= \frac{1}{4}\Delta (P\vert \vert Q)$.

\item[(ii)] \begin{itemize} \item[(a)] $\Omega _0 (P\vert \vert Q)
= F(P\vert \vert Q)$. \item[(b)]$\Omega _0 (Q\vert \vert P) =
F(Q\vert \vert P)$.
\end{itemize}

\item[(iii)] \begin{itemize} \item[(a)] $\Omega _1 (P\vert \vert
Q) = G(P\vert \vert Q)$.

\item[(b)] $\Omega _1 (Q\vert \vert P) = G(Q\vert \vert P)$.
\end{itemize}
\item[(iv)] \begin{itemize}\item[(a)] $\Omega _2 (P\vert \vert Q)
= \frac{1}{8}\chi ^2(Q\vert \vert P)$.

\item[(b)] $\Omega _2 (Q\vert \vert P) = \frac{1}{8}\chi ^2(P\vert
\vert Q)$.
\end{itemize}
\end{itemize}

The expression $\Delta (P\vert \vert Q)$ appearing in part (i) is
the well known \textit{triangular discrimination}, and is given by
\begin{equation}
\label{eq16} \Delta (P\vert \vert Q) = \sum\limits_{i = 1}^n
{\frac{(p_i - q_i )^2}{p_i + q_i }} .
\end{equation}

\bigskip
\textbf{$\bullet$ Relative J-Divergence of Type s}

\bigskip
We now propose the following one parameter generalization of the
\textit{relative J-divergence} measures given by (\ref{eq5}).
\begin{equation}
\label{eq17} \zeta _s (P\vert \vert Q) = \begin{cases}
 {D_s (P\vert \vert Q) = (s - 1)^{ - 1}\sum\limits_{i = 1}^n {(p_i - q_i
)\left( {\frac{p_i + q_i }{2q_i }} \right)^{s - 1},} } & {s \ne 1} \\
 {D(P\vert \vert Q) = \sum\limits_{i = 1}^n {(p_i - q_i )\ln \left(
{\frac{p_i + q_i }{2q_i }} \right),} } & {s = 1} \\
\end{cases}
\end{equation}

The adjoint of $\zeta _s (P\vert \vert Q)$ written as $\zeta _s
(Q\vert \vert P)$ is obtained by interchanging $P$ and $Q$, and
$p_i $ and $q_i $ in the expression (\ref{eq17}).\\

These admit the following particular cases:
\begin{itemize}
\item[(i)] $\zeta _0 (P\vert \vert Q) = \zeta _0 (Q\vert \vert P)
= \Delta (P\vert \vert Q)$.

\item[(ii)] \begin{itemize} \item[(a)] $\zeta _1 (P\vert \vert Q)
= D(P\vert \vert Q)$.

\item[(b)] $\zeta _1 (Q\vert \vert P) = D(Q\vert \vert P)$.
\end{itemize}
\item[(iii)]
\begin{itemize}
\item[(a)] $\zeta _2 (P\vert \vert Q) = \frac{1}{2}\chi^{2}(P\vert
\vert Q)$.

\item[(b)] $\zeta _2 (Q\vert \vert P) = \frac{1}{2}\chi^{2}(Q\vert
\vert P)$.
\end{itemize}
\end{itemize}

\bigskip
We observe that the \textit{relative information of type s}, $\Phi
_s (P\vert \vert Q)$, contains, in particular, the classical
measures Bhattacharyya\textit{ coefficient, }$\chi ^2 -
$\textit{divergence }and Hellingar\textit{ discrimination}. The
\textit{unified relative JS and AG -- divergences of type s},
$\Omega _s (P\vert \vert Q)$ and $\Omega _s (Q\vert \vert P)$,
contains, in particular, \textit{triangular discrimination} and
$\chi ^2 - $\textit{divergence}, while the \textit{relative
J-divergences of type s}, $\zeta _s (P\vert \vert Q)$ and $\zeta
_s (Q\vert \vert P)$, yield, in particular, \textit{triangular
discrimination} and $\chi ^2 - $\textit{divergence}.\\

In this paper our aim is to relate these \textit{generalized
measures of type s} with one another. In order to do so, we make
use of the Csisz\'{a}r \textit{f-divergence} and its properties.

\section{Csisz\'{a}r $f-$Divergence and its
Particular Cases}

For a function $f:[0,\infty ) \to \mathbb{R}$, the
\textit{f-divergence} measure introduced by Csisz\'{a}r
\cite{csi1} is given by
\begin{equation}
\label{eq18}
C_f (P\vert \vert Q) =
\sum\limits_{i = 1}^n {q_i f\left( {\frac{p_i }{q_i }} \right)} ,
\end{equation}

\noindent for all $P,Q \in \Gamma _n $.

The following result is well known in the literature.

\begin{theorem} \label{the31} (Csisz\'{a}r \cite{csi1, csi2}). If the
function $f$ is convex and normalized, i.e., $f(1) = 0$, then $C_f
(P\vert \vert Q)$ and its adjoint $C_f (Q\vert \vert P)$ are both
nonnegative and convex in the pair of probability distribution
$(P,Q) \in \Gamma _n \times \Gamma _n $.
\end{theorem}

The generalized measures given in Section 2 can be written as
particular cases of Csisz\'{a}r \textit{f-divergence}
(\ref{eq18}). These particular cases are given by the following
examples.

\begin{example} \label{exa31} (\textit{Relative information of type
s}). Let us consider
\begin{equation}
\label{eq19} \phi _s (x) = \begin{cases}
 {\left[ {s(s - 1)} \right]^{ - 1}\left[ {x^s - 1 - s(x - 1)} \right],} & {s
\ne 0,1} \\
 {x - 1 - \ln x,} & {s = 0} \\
 {1 - x + x\ln x,} & {s = 1} \\
\end{cases},
\end{equation}

\noindent for all $x > 0$ in (\ref{eq18}). Then $C_f (P\vert \vert
Q) = \Phi _s (P\vert \vert Q)$.
\end{example}

\bigskip
\begin{example} \label{exa32}(\textit{Relative JS and AG -- divergence
of type s}). Let us consider
\begin{equation}
\label{eq20} \psi _s (x) = \begin{cases}
 {\left[ {s(s - 1)} \right]^{ - 1}\left[ {x\left( {\frac{x + 1}{2x}}
\right)^s - x - s\left( {\frac{1 - x}{2}} \right)} \right],} & {s \ne 0,1}
\\
 {\frac{1 - x}{2} - x\ln \left( {\frac{x + 1}{2x}} \right),} & {s = 0} \\
 {\frac{x - 1}{2} + \left( {\frac{x + 1}{2}} \right)\ln \left( {\frac{x +
1}{2x}} \right),} & {s = 1} \\
\end{cases},
\end{equation}

\noindent for all $x > 0$ in (\ref{eq18}). Then $C_f (P\vert \vert
Q) = \Omega _s (P\vert \vert Q)$.
\end{example}

\bigskip
\begin{example} \label{exa33} (\textit{Adjoint of Relative JS and AG
-- divergence of type s}). Let us consider
\begin{equation}
\label{eq21} \upsilon _s (x) = \begin{cases}
 {\left[ {s(s - 1)} \right]^{ - 1}\left[ {\left( {\frac{x + 1}{2}} \right)^s
- 1 - s\left( {\frac{x - 1}{2}} \right)} \right],} & {s \ne 0,1} \\
 {\frac{x - 1}{2} + \ln \left( {\frac{2}{x + 1}} \right),} & {s = 0} \\
 {\frac{1 - x}{2} + \frac{x + 1}{2}\ln \left( {\frac{x + 1}{2}} \right),} &
{s = 1} \\
\end{cases},
\end{equation}

\noindent for all $x > 0$ in (\ref{eq18}). Then $C_f (P\vert \vert
Q) = \Omega _s (Q\vert \vert P)$.
\end{example}

\bigskip
\begin{example} \label{exa34}(\textit{Relative J-divergence of type
s}). Let us consider
\begin{equation}
\label{eq22} \xi _s (x) = \begin{cases}
 {(s - 1)^{ - 1}(x - 1)\left[ {\left( {\frac{x + 1}{2}} \right)^{s - 1} - 1}
\right],} & {s \ne 1} \\
 {(x - 1)\ln \left( {\frac{x + 1}{2}} \right),} & {s = 1} \\
\end{cases},
\end{equation}

\noindent for all $x > 0$ in (\ref{eq18}). Then $C_f (P\vert \vert
Q) = \zeta _s \left( {P\vert \vert Q} \right)$.
\end{example}

\bigskip
\begin{example} \label{exa35} (\textit{Adjoint of relative
J-divergence of type s}). Let us consider
\begin{equation}
\label{eq23} \varsigma _s (x) = \begin{cases}
 {(s - 1)^{ - 1}(1 - x)\left[ {\left( {\frac{x + 1}{2x}} \right)^{s - 1} -
1} \right],} & {s \ne 1} \\
 {(1 - x)\ln \left( {\frac{x + 1}{2x}} \right),} & {s = 1} \\
\end{cases},
\end{equation}

\noindent for all $x > 0$ in (\ref{eq18}). Then $C_f (P\vert \vert
Q) = \zeta _s \left( {Q\vert \vert P} \right)$.
\end{example}

By considering the second order derivative of the functions given
by (\ref{eq19})-(\ref{eq23}) with respect to $x$, and applying the
Theorem \ref{the31}, it can easily be checked that the measures
$\Phi _s (P\vert \vert Q)$, $\Omega _s (P\vert \vert Q)$, $\Omega
_s (Q\vert \vert P)$, $\zeta _s \left( {P\vert \vert Q} \right)$
and $\zeta _s \left( {Q\vert \vert P} \right)$ are
\textit{nonnegative} and \textit{convex} in the pair of
probability distributions $(P,Q) \in \Gamma _n \times \Gamma _n $
respectively, for all $s \in \mathbb{R}$ for the measures $\Phi _s
(P\vert \vert Q)$, $\Omega _s (P\vert \vert Q)$ and $\Omega _s
(Q\vert \vert P)$, and $0 \leqslant s \leqslant 4$ for the
measures $\zeta _s \left( {P\vert \vert Q} \right)$ and $\zeta _s
\left( {Q\vert \vert P} \right)$.

\begin{theorem} \label{the32} (Taneja \cite{tan8}). Let $f_1 ,f_2 :I
\subset \mathbb{R}_ + \to \mathbb{R}$ be two normalized functions,
i.e., $f_1 (1) = f_2 (1) = 0$ and satisfy the assumptions:

(i) $f_1 $ and $f_2 $ are twice differentiable on $(r,R)$;

(ii) there exists the real constants $m,M$ such that $0 \leqslant
m < M$ and
\begin{equation}
\label{eq24}
m \leqslant \frac{f_1 ^{\prime \prime }(x)}{f_2 ^{\prime \prime }(x)}
\leqslant M,
\quad
f_2 ^{\prime \prime }(x) > 0,
\quad
\forall x \in (r,R)
\end{equation}

\noindent then we have the inequalities:
\begin{equation}
\label{eq25}
m\mbox{ }C_{f_2 } (P\vert \vert Q) \leqslant C_{f_1 } (P\vert \vert Q)
\leqslant M\mbox{ }C_{f_2 } (P\vert \vert Q).
\end{equation}
\end{theorem}

\begin{proof} Let us consider the functions $\eta _{m.s} (
\cdot )$ and $\eta _{M.s} ( \cdot )$ given by
\begin{equation}
\label{eq26}
\eta _m (x) = f_1 (x) - m\mbox{ }f_2 (x)
\end{equation}

\noindent and
\begin{equation}
\label{eq27}
\eta _M (x) = M\mbox{ }f_2 (x) - f_1 (x),
\end{equation}

\noindent respectively, where $m$ and $M$ are as given by
(\ref{eq24})

Since $f_1 (x)$ and $f_2 (x)$ are normalized, i.e., $f_1 (1) = f_2
(1) = 0$, then $\eta _m ( \cdot )$ and $\eta _M ( \cdot )$ are
also normalized, i.e., $\eta _m (1) = 0$ and $\eta _M (1) = 0$.
Also, the functions $f_1 (x)$ and $f_2 (x)$ are twice
differentiable. Then in view of (\ref{eq24}), we have
\begin{equation}
\label{eq28}
{\eta }''_m (x) = f_1 ^{\prime \prime }(x) - m\mbox{ }f_2 ^{\prime \prime
}(x)
 = f_2 ^{\prime \prime }(x)\left( {\frac{f_1 ^{\prime \prime }(x)}{f_2
^{\prime \prime }(x)} - m} \right) \geqslant 0
\end{equation}

\noindent and
\begin{equation}
\label{eq29}
{\eta }''_M (x) = M\mbox{ }f_2 ^{\prime \prime }(x) - f_1 ^{\prime \prime
}(x)
 = f_2 ^{\prime \prime }(x)\left( {M - \frac{f_1 ^{\prime \prime }(x)}{f_2
^{\prime \prime }(x)}} \right) \geqslant 0,
\end{equation}

\noindent for all $x \in (r,R)$.

In view of (\ref{eq28}) and (\ref{eq29}), we can say that the
functions $\eta _m ( \cdot )$ and $\eta _M ( \cdot )$ are convex
on $(r,R)$.

According to Theorem \ref{the31}, we have
\begin{equation}
\label{eq30}
C_{\eta _m } (P\vert \vert Q) = C_{f_1 - mf_2 } (P\vert \vert Q) = C_{f_1 }
(P\vert \vert Q) - m\mbox{ }C_{f_2 } (P\vert \vert Q) \geqslant 0,
\end{equation}

\noindent and
\begin{equation}
\label{eq31}
C_{\eta _M } (P\vert \vert Q) = C_{Mf_2 - f_1 } (P\vert \vert Q) = M\mbox{
}C_{f_2 } (P\vert \vert Q) - C_{f_1 } (P\vert \vert Q) \geqslant 0.
\end{equation}

Combining (\ref{eq30}) and (\ref{eq31}), we get (\ref{eq25}).
\end{proof}

\section{Inequalities Among Generalized Relative Divergences}

In this section we shall relate the \textit{relative divergence
measures of type s} given in Section 2. The main results of this
paper are summarized in the following theorem.

\begin{theorem} \label{the41} Let the generating functions given by
(\ref{eq19})-(\ref{eq23}) are twice differentiable in interval
$(r,R)$ with $0 < r \leqslant R$. Then, we have the following
inequalities among the \textit{generalized measures}:
\begin{itemize}
\item[(i)] $\Omega _s (P\vert \vert Q)$ and $\Phi _t (P\vert \vert
Q)$:
\begin{align}
\label{eq32} \frac{1}{4r^{t + 1}}& \left( {\frac{r + 1}{2r}}
\right)^{s - 2}\Phi _t (P\vert \vert Q) \\
& \leqslant \Omega _s (P\vert \vert Q) \leqslant \frac{1}{4R^{t +
1}}\left( {\frac{R + 1}{2R}} \right)^{s - 2}\Phi _t (P\vert \vert
Q), \,\, s + t \leqslant 1, \,\, t \leqslant - 1\notag
\end{align}

\noindent and
\begin{align}
\label{eq33} \frac{1}{4R^{t + 1}}& \left( {\frac{R + 1}{2R}}
\right)^{s - 2}\Phi _t (P\vert \vert Q) \\
&\leqslant \Omega _s (P\vert \vert Q)
 \leqslant \frac{1}{4r^{t + 1}}\left( {\frac{r + 1}{2r}}
\right)^{s - 2}\Phi _t (P\vert \vert Q), \,\, s + t \geqslant 1,
\,\, t \geqslant - 1.\notag
\end{align}

\item[(ii)] $\Omega _s (Q\vert \vert P)$ and $\Phi _t (P\vert
\vert Q)$:
\begin{align}
\label{eq34} \frac{1}{4r^{t - 2}}& \left( {\frac{r + 1}{2}}
\right)^{s - 2}\Phi _t (P\vert \vert Q)\\
 & \leqslant \Omega _s (Q\vert \vert P) \leqslant \frac{1}{4R^{t - 2}}\left( {\frac{R + 1}{2}}
\right)^{s - 2}\Phi _t (P\vert \vert Q), \,\, s \geqslant t,\,t
\leqslant 2\notag
\end{align}

\noindent and
\begin{align}
\label{eq35} \frac{1}{4R^{t - 2}}& \left( {\frac{R + 1}{2}}
\right)^{s - 2}\Phi _t (P\vert \vert Q) \\
& \leqslant \Omega _s (Q\vert \vert P) \leqslant \frac{1}{4r^{t -
2}}\left( {\frac{r + 1}{2}} \right)^{s - 2}\Phi _t (P\vert \vert
Q), \,\, s \leqslant t,\,t \geqslant 2.\notag
\end{align}

\item[(iii)] $\zeta _s (P\vert \vert Q)$ and $\Phi _t (P\vert
\vert Q)$:
\begin{align}
\label{eq36} \left( {\frac{r + 1}{2}} \right)^{s - 3}& \left(
{\frac{sr + 4 - s}{4r^{t - 2}}} \right)\Phi _t (P\vert \vert Q)\\
& \leqslant \zeta _s (P\vert \vert Q) \leqslant \left( {\frac{R +
1}{2}} \right)^{s - 3}\left( {\frac{sR + 4 - s}{4R^{t - 2}}}
\right)\Phi _t (P\vert \vert
Q),\notag \\
& \qquad \qquad \qquad \qquad \qquad 0 \leqslant s \leqslant 4,\,t
\leqslant 2,\,s \geqslant t + 1\notag
\end{align}

\noindent and
\begin{align}
\label{eq37} \left( {\frac{R + 1}{2}} \right)^{s - 3} &
\left(\frac{sR + 4 - s}{4R^{t - 2}}\right) \Phi _t (P\vert \vert
Q) \\
& \leqslant \zeta _s (P\vert \vert Q) \leqslant \left( {\frac{r +
1}{2}} \right)^{s - 3} \left( \frac{sr + 4 - s}{4r^{t - 2}}\right)
\Phi _t (P\vert \vert
Q),\notag\\
& \qquad \qquad \qquad \qquad \qquad 0 \leqslant s \leqslant 4,\,t
\geqslant 2,\,s \leqslant t + 1.\notag
\end{align}

\item[(iv)] $\zeta _s (Q\vert \vert P)$ and $\Phi _t (P\vert \vert
Q)$:
\begin{align}
\label{eq38} \left( {\frac{r + 1}{2r}} \right)^{s - 3}&
\left(\frac{(4 - s)r + s}{4r^{t + 2}}\right)\Phi _t (P\vert \vert
Q) \\
& \leqslant \zeta _s (Q\vert \vert P) \leqslant \left( {\frac{R +
1}{2R}} \right)^{s - 3}\left(\frac{(4 - s)R + s}{4R^{t +
2}}\right)\Phi _t (P\vert \vert Q), \notag\\
& \qquad \qquad \qquad \qquad \qquad 0 \leqslant s \leqslant 4,
\,\, t \leqslant - 1, \,\, s + t \leqslant 1\notag
\end{align}

\noindent and
\begin{align}
\label{eq39} \left( {\frac{R + 1}{2R}} \right)^{s - 3}&
\left(\frac{(4 - s)R + s}{4R^{t + 2}}\right)\Phi _t (P\vert \vert
Q)\\
& \leqslant \zeta _s (Q\vert \vert P) \leqslant \left( {\frac{r +
1}{2r}} \right)^{s - 3}\left(\frac{(4 - s)r + s}{4r^{t +
2}}\right)\Phi _t (P\vert \vert Q),\notag\\
& \qquad \qquad \qquad \qquad \qquad 0 \leqslant s \leqslant 4,
\,\, t \geqslant - 1, \,\, s + t \geqslant 2.\notag
\end{align}

\item[(v)] $\Omega _s (Q\vert \vert P)$ and $\Omega _t (P\vert
\vert Q)$:
\begin{align}
\label{eq40} r^{t + 1}& \left( {\frac{r + 1}{2}} \right)^{s -
t}\Omega _t (P\vert \vert Q) \\
& \leqslant \Omega _s (Q\vert \vert P) \leqslant R^{t + 1}\left(
{\frac{R + 1}{2}} \right)^{s - t}\Omega _t (P\vert \vert Q), \,\,
s \geqslant - 1, \,\, t \geqslant - 1.\notag
\end{align}

\noindent and
\begin{align}
\label{eq41} R^{t + 1}& \left( {\frac{R + 1}{2}} \right)^{s -
t}\Omega _t (P\vert \vert Q) \\
&  \leqslant \Omega _s (Q\vert \vert P) \leqslant r^{t + 1}\left(
{\frac{r + 1}{2}} \right)^{s - t}\Omega _t (P\vert \vert Q), \,\,
s \leqslant - 1, \,\, t \leqslant - 1.\notag
\end{align}

\item[(vi)] $\zeta _s (P\vert \vert Q)$ and $\Omega _t (P\vert
\vert Q)$:
\begin{align}
\label{eq42} & r^{t + 1}\left( {\frac{r + 1}{2}} \right)^{s - t -
1} \left( {sr + 4 - s} \right)\Omega _t (P\vert \vert Q)\\
& \leqslant \zeta _s (P\vert \vert Q)\leqslant R^{t + 1}\left(
{\frac{R + 1}{2}} \right)^{s - t - 1}\left( {sR + 4 - s}
\right)\Omega _t (P\vert \vert Q), \,\, 0 \leqslant s \leqslant 4,
\,\, t \geqslant - 1.\notag
\end{align}

\item[(vii)] $\zeta _s (Q\vert \vert P)$ and $\Omega _t (P\vert
\vert Q)$:
\begin{align}
\label{eq43} \left( {\frac{r + 1}{2r}} \right)^{s - t - 1}& \left(
{\frac{(4 - s)r + s}{r}} \right)\Omega _t (P\vert \vert Q)\\
& \leqslant \zeta _s (Q\vert \vert P) \leqslant \left( {\frac{R +
1}{2R}} \right)^{s - t - 1}\left( {\frac{(4 - s)R + s}{R}}
\right)\Omega _t (P\vert \vert Q), \notag\\
& \qquad \qquad \qquad \qquad 0 \leqslant s \leqslant 4,\,\,t
\geqslant s,\,\,t(4 - s) \geqslant 6s - s^2 - 4 \notag
\end{align}

\noindent and
\begin{align}
\label{eq44} \left( {\frac{R + 1}{2R}} \right)^{s - t - 1}& \left(
{\frac{(4 - s)R + s}{R}} \right)\Omega _t (P\vert \vert Q)\\
& \leqslant \zeta _s (Q\vert \vert P)  \leqslant \left( {\frac{r +
1}{2r}} \right)^{s - t - 1}\left( {\frac{(4 - s)r + s}{r}}
\right)\Omega _t (P\vert \vert Q),\notag\\
& \qquad \qquad \qquad \qquad 0 \leqslant s \leqslant 4,\,\,t
\leqslant s,\,\,t(4 - s) \leqslant 6s - s^2 - 4 \notag
\end{align}

\item[(viii)] $\zeta _s (P\vert \vert Q)$ and $\Omega _t (Q\vert
\vert P)$:
\begin{align}
\label{eq45} \left( {\frac{r + 1}{2}} \right)^{s - t - 1}& \left(
{sr + 4 - s} \right)\Omega _t (Q\vert \vert P) \\
& \leqslant \zeta _s (P\vert \vert Q)\leqslant \left( {\frac{R +
1}{2}} \right)^{s - t - 1}\left[ {sR + (4 - s)} \right]\Omega _t
(Q\vert \vert P),\notag\\
& \qquad \qquad \qquad \qquad 0 \leqslant s \leqslant 4,\,\,s
\geqslant t,\,\,s(t - s + 6) \geqslant 4(1 + t)\notag
\end{align}

\noindent and
\begin{align}
\label{eq46} \left( {\frac{R + 1}{2}} \right)^{s - t - 1}& \left(
{sR + 4 - s} \right)\Omega _t (Q\vert \vert P) \\
& \leqslant \zeta _s (P\vert \vert Q) \leqslant \left( {\frac{r +
1}{2}} \right)^{s - t - 1}\left( {sr + 4 - s} \right)\Omega _t
(Q\vert \vert P),\notag \\
& \qquad \qquad \qquad \qquad 0 \leqslant s \leqslant 4,\,\,s
\leqslant t,\,\,s(t - s + 6) \leqslant 4(1 + t).\notag
\end{align}

\item[(ix)] $\zeta _s (Q\vert \vert P)$ and $\Omega _t (Q\vert
\vert P)$:
\begin{align}
\label{eq47} & \frac{1}{R^{s + 1}} \left( {\frac{R + 1}{2}}
\right)^{s - t - 1}\left[ {(4 - s)R + s} \right]\Omega _t (Q\vert
\vert P) \\
& \leqslant \zeta _s (Q\vert \vert P) \leqslant \frac{1}{r^{s +
1}}\left( {\frac{r + 1}{2}} \right)^{s - t - 1}\left[ {(4 - s)R +
s} \right]\Omega _t (Q\vert \vert P), \,0 \leqslant s \leqslant 4,
\, t \geqslant - 1.\notag
\end{align}

\item[(x)] $\zeta _s (Q\vert \vert P)$ and $\zeta _t (P\vert \vert
Q)$:
\begin{align}
\label{eq48} & \frac{1}{R^{s + 1}}\left( {\frac{R + 1}{2}}
\right)^{s - t}\left( {\frac{(4 - s)R + s}{tR + 4 - t}}
\right)\zeta _t (P\vert \vert Q) \\
& \leqslant \zeta _s (Q\vert \vert P) \leqslant \frac{1}{r^{s +
1}}\left( {\frac{r + 1}{2}} \right)^{s - t}\left( {\frac{(4 - s)r
+ s}{tr + 4 - t}} \right)\zeta _t (P\vert \vert Q), \,\, 2
\leqslant s \leqslant 4,\,\,2 \leqslant t \leqslant 4.\notag
\end{align}
\end{itemize}
\end{theorem}

\begin{proof} (i) Let us consider
\begin{equation}
\label{eq49}
g_{(\psi _s ,\phi _t )} (x) = \frac{{\psi }''_s (x)}{{\phi }''_t (x)} =
\frac{\frac{1}{4x^3}\left( {\frac{x + 1}{2x}} \right)^{s - 2}}{x^{t - 2}}
 = \frac{1}{4x^{t + 1}}\left( {\frac{x + 1}{2x}} \right)^{s - 2},
\end{equation}

\noindent for all $x \in (0,\infty )$.

From (\ref{eq49}) one has
\begin{equation}
\label{eq50}
{g}'_{(\psi _s ,\phi _t )} (x) = - \left( {\frac{x + 1}{2x}} \right)^{s -
2}\frac{x(t + 1) + t + s - 1}{4x^{t + 2}(x + 1)}
\begin{cases}
 { \geqslant 0,} & {t \leqslant - 1,\mbox{ }s + t \leqslant 1} \\
 { \leqslant 0,} & {t \geqslant - 1,\mbox{ }s + t \geqslant 1} \\
\end{cases}.
\end{equation}

In view of (\ref{eq50}) we conclude the followings
\begin{equation}
\label{eq51} m = \mathop {\inf }\limits_{x \in [r,R]} g_{(\psi _s
,\phi _t )} (x) = \begin{cases}
 {\frac{1}{4r^{t + 1}}\left( {\frac{r + 1}{2r}} \right)^{s - 2},} & {s + t
\leqslant 1,\mbox{ }t \leqslant - 1} \\
 {\frac{1}{4R^{t + 1}}\left( {\frac{R + 1}{2R}} \right)^{s - 2},} & {s + t
\geqslant 1,\mbox{ }t \geqslant - 1} \\
\end{cases}
\end{equation}

\noindent and
\begin{equation}
\label{eq52} M = \mathop {\sup }\limits_{x \in [r,R]} g_{(\psi _s
,\phi _t )} (x) = \begin{cases}
 {\frac{1}{4R^{t + 1}}\left( {\frac{R + 1}{2R}} \right)^{s - 2},} & {s + t
\leqslant 1,\mbox{ }t \leqslant - 1} \\
 {\frac{1}{4r^{t + 1}}\left( {\frac{r + 1}{2r}} \right)^{s - 2},} & {s + t
\geqslant 1,\mbox{ }t \geqslant - 1} \\
\end{cases}.
\end{equation}

Now (\ref{eq51}) and (\ref{eq52}) together with (\ref{eq25}) give
the inequalities (\ref{eq32}) and (\ref{eq33}).
\end{proof}

The proof of other parts (ii)-(x) follows on similar lines.

\subsection*{4.1. Particular Cases} Here below we have considered
some particular cases of the inequalities
(\ref{eq32})-(\ref{eq48}).

\begin{itemize}
\item Take $t = \frac{1}{2}, s = 2$ in (\ref{eq34}) or in
(\ref{eq36}), one gets
\begin{equation*}
\label{eq53} r \leqslant \frac{\sqrt[3]{\chi ^2(P\vert \vert
Q)^2}}{4 \sqrt[3]{h(P\vert \vert Q)^2}} \leqslant R.
\end{equation*}

\item Take $t = \frac{1}{2}, s = 2$ in (\ref{eq33}) or $t =
\frac{1}{2}, s = 2$ in (\ref{eq39}), one gets
\begin{equation*}
\label{eq54} r \leqslant \frac{4 \sqrt[3]{h(P\vert \vert
Q)^2}}{\sqrt[3]{\chi ^2(Q\vert \vert P)^2}} \leqslant R.
\end{equation*}

\item Take $t = 2, s = 2$ in (\ref{eq33}) or in (\ref{eq39}) or in
(\ref{eq40}) or in (\ref{eq42}) or in (\ref{eq47}) or in
(\ref{eq48}) or $t = - 1,\mbox{ }s = 2$ in (\ref{eq34}) or in
(\ref{eq36}), one gets
\begin{equation*}
\label{eq55}
r \leqslant \frac{\sqrt[3]{\chi ^2(P\vert \vert Q)}}{\sqrt[3]{\chi ^2(Q\vert
\vert P)}} \leqslant R.
\end{equation*}

\item Take $t = 1, s = 2$ in (\ref{eq34}) or in (\ref{eq36}), one
gets
\begin{equation*}
\label{eq56}
r \leqslant \frac{\chi ^2(P\vert \vert Q)}{2K(P\vert \vert Q)} \leqslant R.
\end{equation*}

\item Take $t = 0, s = 2$ in (\ref{eq33}) or in (\ref{eq39}), one
gets
\begin{equation*}
\label{eq57}
r \leqslant \frac{2K(Q\vert \vert P)}{\chi ^2(Q\vert \vert P)} \leqslant R.
\end{equation*}

\item Take $t = 1, s = 2$ in (\ref{eq33}) or in (\ref{eq39}), one
gets
\begin{equation*}
\label{eq58}
r \leqslant \frac{\sqrt {2K(P\vert \vert Q)} }{\sqrt {\chi ^2(Q\vert \vert
P)} } \leqslant R.
\end{equation*}

\item Take $t = 0, s = 2$ in (\ref{eq34}) or in (\ref{eq36}), one
gets
\begin{equation*}
\label{eq59}
r \leqslant \frac{\sqrt {\chi ^2(P\vert \vert Q)} }{\sqrt {2K(Q\vert \vert
P)} } \leqslant R.
\end{equation*}

\item Take $t = 0, s = 0$ in (\ref{eq40}), one gets
\begin{equation*}
\label{eq60}
r \leqslant \frac{F(Q\vert \vert P)}{F(P\vert \vert Q)} \leqslant R.
\end{equation*}

\item Take $t = 1, s = 1$ in (\ref{eq40}), one gets
\begin{equation*}
\label{eq61}
r \leqslant \frac{\sqrt {G(Q\vert \vert P)} }{\sqrt {G(P\vert \vert Q)} }
\leqslant R.
\end{equation*}

\item Take $t = 2, s = - 1$ in (\ref{eq33}) or $t = 2,\mbox{ }s =
0$ in (\ref{eq39}) or (\ref{eq46}) or in (\ref{eq47}) or $t = -
1,\mbox{ }s = 2$ in (\ref{eq40}) or in (\ref{eq42}), one gets
\begin{equation*}
\label{eq62}
r \leqslant \frac{\sqrt[3]{4\chi ^2(P\vert \vert Q)} - \sqrt[3]{\Delta
(P\vert \vert Q)}}{\sqrt[3]{\Delta (P\vert \vert Q)}} \leqslant R.
\end{equation*}

\item Take $t = - 1, s = - 1$ in (\ref{eq32}) or in (\ref{eq34})
$t = - 1, s = 0$ in (\ref{eq36}) or $t = 2,\mbox{ }s = - 1$ in
(\ref{eq40}) or $t = 2,\mbox{ }s = 0$ in (\ref{eq42}) or in
(\ref{eq44}) or $t = - 1,\mbox{ }s = 2$ in (\ref{eq44}) or in
(\ref{eq47}), one gets
\begin{equation*}
\label{eq63}
r \leqslant \frac{\sqrt[3]{\Delta (P\vert \vert Q)}}{\sqrt[3]{4\chi
^2(Q\vert \vert P)} - \sqrt[3]{\Delta (P\vert \vert Q)}} \leqslant R.
\end{equation*}

\item Take $t = 0,\mbox{ }s = 1$ in (\ref{eq33}), one gets
\begin{equation*}
\label{eq64}
r \leqslant \frac{K(Q\vert \vert P) - 2G(P\vert \vert Q)}{2G(P\vert \vert
Q)} \leqslant R.
\end{equation*}

\item Take $t = 1, s = 1$ in (\ref{eq34}), one gets
\begin{equation*}
\label{eq65}
r \leqslant \frac{2G(Q\vert \vert P)}{K(P\vert \vert Q) - 2G(Q\vert \vert
P)} \leqslant R.
\end{equation*}

\item Take $t = 1, s = 0$ in (\ref{eq33}), one gets
\begin{equation*}
\label{eq66}
r \leqslant \frac{\sqrt {K(P\vert \vert Q)} - \sqrt {F(P\vert \vert Q)}
}{\sqrt {F(P\vert \vert Q)} } \leqslant R.
\end{equation*}

\item Take $t = 0, s = 0$ in (\ref{eq34}), one gets
\begin{equation*}
\label{eq67}
r \leqslant \frac{\sqrt {F(Q\vert \vert P)} }{\sqrt {K(Q\vert \vert P)} -
\sqrt {F(Q\vert \vert P)} } \leqslant R.
\end{equation*}

\item Take $t = 1, s = - 1$ in (\ref{eq40}) or in $t = 1,\mbox{ }s
= 0$ in (\ref{eq42}) or (\ref{eq43}), one gets
\begin{equation*}
\label{eq68}
r \leqslant \frac{\sqrt {\Delta (P\vert \vert Q)} }{4\sqrt {G(P\vert \vert
Q)} - \sqrt {\Delta (P\vert \vert Q)} } \leqslant R.
\end{equation*}

\item Take $t = - 1, s = 1$ in (\ref{eq40}) or $t = 1, s = 0$ in
(\ref{eq46}) or in (\ref{eq47}), one gets
\begin{equation*}
\label{eq69}
r \leqslant \frac{4\sqrt {G(Q\vert \vert P)} - \sqrt {\Delta (P\vert \vert
Q)} }{\sqrt {\Delta (P\vert \vert Q)} } \leqslant R.
\end{equation*}

\item Take $t = 0, s = - 1$ in (\ref{eq40}) or $t = 0, s = 0$ in
(\ref{eq42}) or in (\ref{eq44}), one gets
\begin{equation*}
\label{eq70}
r \leqslant \frac{\Delta (P\vert \vert Q)}{8F(P\vert \vert Q) - \Delta
(P\vert \vert Q)} \leqslant R.
\end{equation*}

\item Take $t = - 1, s = 0$ in (\ref{eq40}) or $t = 0, s = 0$ in
(\ref{eq46}) or in (\ref{eq47}), one gets
\begin{equation*}
\label{eq71}
r \leqslant \frac{8F(Q\vert \vert P) - \Delta (P\vert \vert Q)}{\Delta
(P\vert \vert Q)} \leqslant R.
\end{equation*}

\item Take $t = 1, s = 1$ in (\ref{eq46}) , one gets
\begin{equation*}
\label{eq75} r \leqslant \frac{6G(Q\vert \vert P) - D(P\vert \vert
Q)}{D(P\vert \vert Q) - 2G(Q\vert \vert P)} \leqslant R.
\end{equation*}

\item Take $t = 1, s = 1$ in (\ref{eq44}) , one gets
\begin{equation*}
\label{eq76} r \leqslant \frac{D(Q\vert \vert P) - 2\,G(P\vert
\vert Q)}{6\,G(P\vert \vert Q) - D(Q\vert \vert P)} \leqslant R.
\end{equation*}

\item Take $t = - 1, s = 1$ in (\ref{eq32}) or $t = 1, s = 2$ in
(\ref{eq44}), one gets
\begin{equation*}
\label{eq74} r \leqslant \frac{4G(P\vert \vert Q)}{\chi ^2(Q\vert
\vert P) - 4G(P\vert \vert Q)} \leqslant R.
\end{equation*}

\item Take $t = 0, s = 1$ in (\ref{eq32}) or $t = 1, s = 2$ in
(\ref{eq44}), one gets
\begin{equation*}
\label{eq74} r \leqslant \frac{F(P\vert \vert Q)}{D(Q\vert \vert
P) - 3F(P\vert \vert Q)} \leqslant R.
\end{equation*}

\item Take $t = - 1, s = 0$ in (\ref{eq32}) or $t = 0, s = 2$ in
(\ref{eq44}), one gets
\begin{equation*}
\label{eq72} r \leqslant \frac{\sqrt {2F(P\vert \vert Q)} }{\sqrt
{\chi ^2(Q\vert \vert P)} - \sqrt {2F(P\vert \vert Q)} } \leqslant
R.
\end{equation*}

\item Take $t = 0, s = 1$ in (\ref{eq42}) , one gets
\begin{equation*}
\label{eq79} r \leqslant \frac{\sqrt {4D(P\vert \vert Q) +
9F(P\vert \vert Q} - 3\sqrt {F(P\vert \vert Q)} }{2\sqrt {F(P\vert
\vert Q)} } \leqslant R.
\end{equation*}

\item Take $t = 0,s = 1$ in (\ref{eq47}) , one gets
\begin{equation*}
\label{eq81} r \leqslant \frac{2\sqrt {F(Q\vert \vert P)} }{\sqrt
{4D(Q\vert \vert P) + 9F(Q\vert \vert P)} - 3\sqrt {F(Q\vert \vert
P)} } \leqslant R.
\end{equation*}

\item Take $t = 1, s = 0$ in (\ref{eq40}) , one gets
\begin{equation*}
\label{eq77} r \leqslant \frac{2\sqrt {F(Q\vert \vert P)} }{\sqrt
{8G(P\vert \vert Q) + F(Q\vert \vert P)} - \sqrt {F(Q\vert \vert
P)} } \leqslant R.
\end{equation*}

\item Take $t = 0, s = 1$ in (\ref{eq40}), one gets
\begin{equation*}
\label{eq80} r \leqslant \frac{\sqrt {8G(Q\vert \vert P) +
F(P\vert \vert Q)} - \sqrt {F(P\vert \vert Q)} }{2\sqrt {F(P\vert
\vert Q)} } \leqslant R.
\end{equation*}

\item Take $t = 0, s = 1$ in (\ref{eq34}) , one gets
\begin{equation*}
\label{eq78} r \leqslant \frac{2\sqrt {G(Q\vert \vert P)} }{\sqrt
{2K(Q\vert \vert P) + G(Q\vert \vert P)} - \sqrt {G(Q\vert \vert
P)} } \leqslant R.
\end{equation*}

\item Take $t = 1, s = 1$ in (\ref{eq33}), one gets
\begin{equation*}
\label{eq83} r \leqslant \frac{\sqrt {2K(P\vert \vert Q) +
G(P\vert \vert Q)} - \sqrt {G(P\vert \vert Q)} }{2\sqrt {G(P\vert
\vert Q)} } \leqslant R.
\end{equation*}

\item Take $t = - 1, s = 1$ in (\ref{eq42}), one gets
\begin{equation*}
\label{eq84} r \leqslant \frac{\sqrt {8D(P\vert \vert Q) + \Delta
(P\vert \vert Q)} - 2\sqrt {\Delta (P\vert \vert Q)} }{\sqrt
{\Delta (P\vert \vert Q)} } \leqslant R.
\end{equation*}

\item Take $t = - 1, s = 1$ in (\ref{eq44}) or in (\ref{eq47}),
one gets
\begin{equation*}
\label{eq85} r \leqslant \frac{\sqrt {\Delta (P\vert \vert Q)}
}{\sqrt {8D(Q\vert \vert P) + \Delta (P\vert \vert Q)} - 2\sqrt
{\Delta (P\vert \vert Q)} } \leqslant R.
\end{equation*}

\item Take $t = 2, s = 1$ in (\ref{eq46}), one gets
\begin{equation*}
\label{eq82} r \leqslant \frac{5\sqrt {\chi ^2(P\vert \vert Q)} -
\sqrt {16D(P\vert \vert Q) + \chi ^2(P\vert \vert Q)} }{\sqrt
{16D(P\vert \vert Q) + \chi ^2(P\vert \vert Q)} - \sqrt {\chi
^2(P\vert \vert Q)} } \leqslant R.
\end{equation*}
\end{itemize}

\bigskip
\begin{remark} \label{rem41} \begin{itemize} \item[(i)]
The inequalities (\ref{eq32})-(\ref{eq48}) admits much more
particular cases but we have specified here only the one that can
be written in simplified form.

\item[(ii)] There are some other similar kind of relations that
can't be obtained from the inequalities (\ref{eq32})-(\ref{eq48})
such as: between $K(P\vert \vert Q)$ and $K(Q\vert \vert P)$,
between $F(P\vert \vert Q)$ and $G(P\vert \vert Q)$, and between
$F(Q\vert \vert P)$ and $G(Q\vert \vert P)$, between $D(P\vert
\vert Q)$ and $D(Q\vert \vert P)$, etc. These can be seen in
Taneja \cite{tan7, tan9} and Taneja and Kumar \cite{tak}.
\end{itemize}
\end{remark}

\bigskip
\begin{center}
Acknowledgements
\end{center}

\bigskip
The authors are thankful to the referee for valuable comments and
suggestions on an earlier version of the paper.


\begin{thebibliography}{99}
\bibitem{pea} K. PEARSON, On
the Criterion that a given system of deviations from the probable
in the case of correlated system of variables is such that it can
be reasonable supposed to have arisen from random sampling,
\textit{Phil. Mag.}, \textbf{50}(1900), 157-172.
\end{thebibliography}
\end{document}